\documentclass[reqno]{amsart}

\usepackage{hyperref}

\theoremstyle{exercise}

\theoremstyle{definition}

\usepackage{graphicx}
\usepackage{pstricks, enumerate, pst-node, pst-text, pst-plot}
\theoremstyle{remark}

\numberwithin{equation}{section}

\newcommand{\intav}[1]{\mathchoice {\mathop{\vrule width 6pt height 3 pt depth  -2.5pt
\kern -8pt \intop}\nolimits_{\kern -6pt#1}} {\mathop{\vrule width
5pt height 3  pt depth -2.6pt \kern -6pt \intop}\nolimits_{#1}}
{\mathop{\vrule width 5pt height 3 pt depth -2.6pt \kern -6pt
\intop}\nolimits_{#1}} {\mathop{\vrule width 5pt height 3 pt depth
-2.6pt \kern -6pt \intop}\nolimits_{#1}}}

\newcommand{\intavl}[1]{\mathchoice {\mathop{\vrule width 6pt height 3 pt depth  -2.5pt
\kern -8pt \intop}\limits_{\kern -6pt#1}} {\mathop{\vrule width 5pt
height 3  pt depth -2.6pt \kern -6pt \intop}\nolimits_{#1}}
{\mathop{\vrule width 5pt height 3 pt depth -2.6pt \kern -6pt
\intop}\nolimits_{#1}} {\mathop{\vrule width 5pt height 3 pt depth
-2.6pt \kern -6pt \intop}\nolimits_{#1}}}

\begin{document}

\title[Classical spectra and dynamical Cantor sets]{Diophantine approximations, Lagrange and Markov spectra, and dynamical Cantor sets}

\author[Carlos Matheus]{Carlos Matheus}
\address{Centre de Math\'ematiques Laurent Schwartz, CNRS (UMR 7640), \'Ecole Polytechnique, 91128 Palaiseau, France.}
\email{carlos.matheus@math.cnrs.fr}
\thanks{We are thankful to Erica Flapan and William McCallum for their encouragements during the preparation of this text.}
\author[Carlos Gustavo Moreira]{Carlos Gustavo Moreira}
\address{IMPA, Estrada D. Castorina, 110, 22460-320, Rio de Janeiro, Brazil.}
\email{gugu@impa.br}


\date{\today}



\begin{abstract} 
This text is a slightly expanded version of a survey article on certain aspects of low dimensional dynamics and number theory written after a kind invitation by the editors of the Notices of the American Mathematical Society.  
\end{abstract}

\maketitle


\section{Diophantine approximations} The seminal works of Diophantus of Alexandria (circa AD 250) on rational approximations to the solutions of certain algebraic equations began the important subfield of Number Theory called \emph{Diophantine approximation}. Among the basic problems in this topic, one has the question of finding rational numbers $p/q\in\mathbb{Q}$ approximating a given real number $\alpha\in\mathbb{R}$ in such a way that the denominator $q$ is not ``big'' \emph{and} the error $|\alpha-p/q|$ is ``small''. 

\subsection{Rational approximations of $\pi$}

The first few decimal digits of the number $\pi$ are well-known: $\pi=3.1415926\dots$. By definition, this provides some rational approximations of $\pi$ like $314/100$ and $3141592/10^6$. Nonetheless, these fractions are certainly \emph{not} answers to the Diophantine problem posed above because we can get \emph{better} approximations with \emph{smaller} denominators: for instance, Archimedes (circa 250 BC) knew that 
$$\left|\pi-\frac{22}{7}\right|< \frac{1}{700} <\left|\pi-\frac{314}{100}\right|$$ 
and it is possible to check that 
$$\left|\pi-\frac{355}{113}\right|< \frac{1}{3\cdot 10^6}<\left|\pi-\frac{3141592}{10^6}\right|.$$

\subsection{Dirichlet's pigeonhole principle}
The example of the number $\pi$ makes us wonder how small $|\alpha-p/q|$ can be when the denominator $q$ varies in a fixed range $1\leq q\leq Q$. A preliminary answer comes from the following elementary remark. Recall that any real number $\alpha$ lies between two consecutive integers, namely $\lfloor\alpha\rfloor\leq \alpha < \lfloor\alpha\rfloor+1$ where $\lfloor\alpha\rfloor\in\mathbb{Z}$ is the \emph{integer part} of $\alpha$. Therefore, given $\alpha\in\mathbb{R}$ and $q\in\mathbb{N}$, we can find $p\in\mathbb{Z}$ such that $|q\alpha-p|\leq 1/2$, i.e., 
$$\left|\alpha-\frac{p}{q}\right| < \frac{1}{2q}.$$

In 1841, Dirichlet used his famous \emph{pigeonhole principle} to significantly improve upon the elementary statement in the previous paragraph: more concretely, for any irrational number $\alpha$, one has that 
$$\#\left\{\frac{p}{q}\in\mathbb{Q}: \left|\alpha-\frac{p}{q}\right|<\frac{1}{q^2}\right\} = \infty.$$
Indeed, given $Q\in\mathbb{N}$, Dirichlet considered how the $Q+1$ numbers $\{i\alpha\}:=i\alpha-\lfloor i\alpha\rfloor\in [0,1)$, $i=0, \dots, Q$, are distributed across the elements of the partition $[0,1)=\bigcup\limits_{n=1}^Q [(n-1)/Q, n/Q)$ into $Q$ intervals. By the pigeonhole principle, two fractional parts, say $\{i\alpha\}, \{j\alpha\}$, $0\leq i<j\leq Q$, must lie in the same interval, say $[(n-1)/Q, n/Q)$, so that $|\{j\alpha\}-\{i\alpha\}|<1/Q$ and, \emph{a fortiori}, there exists $p\in\mathbb{Z}$ such that 
$$\left|\alpha-\frac{p}{j-i}\right|<\frac{1}{(j-i)Q}\leq \frac{1}{(j-i)^2}.$$

\subsection{Hurwitz theorem}
In 1891, Hurwitz proved\footnote{Actually, this result was first established by Korkine and Zolotarev in 1873.} that Dirichlet's theorem is essentially optimal as far as \emph{all} irrational numbers are concerned: one has   
$$\#\left\{\frac{p}{q}\in\mathbb{Q}: \left|\alpha-\frac{p}{q}\right| < \frac{1}{\sqrt{5}\cdot q^2}\right\}=\infty$$  
for all irrational number $\alpha$, and\footnote{Let $\phi:=(1+\sqrt{5})/2$. If $p/q\in\mathbb{Q}$ and $|q\phi-p|\leq 1$, then $1\leq |q^2+pq-p^2| = |q\phi-p|\cdot |q\phi^{-1}+p|$ and $|q\phi^{-1}+p| = |q(\phi^{-1}+\phi)-q\phi+p|\leq q\sqrt{5}+1$. Thus, $\frac{1}{|q(q\phi-p)|}\leq \sqrt{5}+\frac{1}{q}$.}  
$$\#\left\{\frac{p}{q}\in\mathbb{Q}: \left|\frac{1+\sqrt{5}}{2}-\frac{p}{q}\right| < \frac{1}{(\sqrt{5}+\varepsilon)\cdot q^2}\right\}<\infty$$ 
for all $\varepsilon > 0$.  

\section{Classical spectra}\label{s.classical-spectra} Despite the almost optimality of the Dirichlet theorem, we can ask whether it can be improved for \emph{individual} irrational numbers $\alpha$ by inquiring about the nature of the \emph{best} constant $\ell(\alpha)$ among the quantities $c$ such that 
$$\#\left\{\frac{p}{q}\in\mathbb{Q}: \left|\alpha-\frac{p}{q}\right|<\frac{1}{c\cdot q^2}\right\}=\infty,$$
i.e., $\ell(\alpha):=\limsup\limits_{p,q\to\infty}\frac{1}{q^2|\alpha-p/q|} = \limsup\limits_{p,q\to\infty}\frac{1}{|q(q\alpha-p)|}$.

The \emph{Lagrange spectrum} $L$ is the collection of \emph{finite}\footnote{It is possible to show that $\ell(\alpha)=\infty$ for Lebesgue almost every $\alpha$. Hence, the Lagrange spectrum tries to encode Diophantine properties of irrational numbers beyond the probabilistic dominant regime.} best constants of Diophatine approximation, i.e., 
$$L:=\{\ell(\alpha)<\infty: \alpha\in\mathbb{R}\setminus\mathbb{Q}\}.$$
In this setting, Hurwitz theorem says that the minimum of $L$ is $\min L=\sqrt{5}$. 

\subsection{Beginning of the classical spectra} The Lagrange spectrum was systematically studied\footnote{This work was part of Markov's M.Sc. dissertation (directly motivated by the results of Korkine and Zolotarev), while Markov's PhD thesis focused on the theory of Markov chains in Probability Theory.} in connection with the theory of binary quadratic forms by Markov in 1879. In fact, the quantity $q(q\alpha-p)$ is the value of the binary quadratic form $h_{\alpha}(x,y)=\alpha y^2 - xy$ at the integral point $(p,q)\in\mathbb{Z}^2$, so that the Lagrange spectrum is somewhat related to the \emph{Markov spectrum} $M$ of \emph{finite} best constants 
$$m(h):=\sup\limits_{(p,q)\in\mathbb{Z}^2\setminus\{(0,0)\}} \frac{\sqrt{\Delta(h)}}{|h(p,q)|}$$ 
of Diophantine approximations of real, indefinite, binary\footnote{The theory of $n$-ary quadratic forms is also very interesting, but we will skip it here because the dynamical techniques involved in Margulis' solution of Oppenheim conjecture are dramatically different from the methods appearing in this text.} quadratic forms $h(x,y)=ax^2+bxy+cy^2$ with positive discriminant $\Delta(h)=b^2-4ac>0$. In this context, Markov proved that 
$$L\cap [\sqrt{5},3) = M\cap [\sqrt{5},3) = \left\{\sqrt{9-\frac{4}{z_n^2}}: n\in\mathbb{N}\right\},$$ 
where $z_n$ is a \emph{Markov number}, i.e., the largest coordinate of a solution $(x_n, y_n, z_n)\in\mathbb{N}^3$ of the Markov--Hurwitz equation 
$$x_n^2+y_n^2+z_n^2=3 x_n y_n z_n.$$  

\subsubsection{Fermat's descent on Markov's cubic} The Markov--Hurwitz equation determines a cubic surface $S$ whose integral points are called \emph{Markov triples}. Since the Markov--Hurwitz equation is quadratic on a given variable (when we freeze the other two variables), the cubic surface $S$ has a rich group of automorphisms: besides permuting the coordinates, we can replace $(x,y,z)$ by $(3yz-x,y,z)$, $(x,3xz-y,z)$ or $(x, y, 3xy-z)$ without leaving $S$. The last three automorphisms are called \emph{Vieta involutions} and they were used by Markov to produce a \emph{descent argument} showing that \emph{any} Markov triple $(x,y,z)\in\mathbb{N}^3$ can be obtained from the fundamental solution $(1,1,1)$ after applying a sequence of permutations of coordinates and Vieta involutions. 

In fact, it is not hard to see that a Markov triple $(x,y,z)\in\mathbb{N}^3$ with $x\leq y\leq z$ falls into two categories: either $x=y$ or $x<y<z$. In the first case, $z^2=(3z-2)x^2$, so that $x=y$ divides $z$, say $z= n x$ with $n^2=3nx-2$; this means that $n$ divides $2$, and, thus, $n=1$ and $(x,y,z)=(1,1,1)$ or $n=2$ and $(x,y,z)=(1,1,2)$. In the second case, we apply Vieta involution $(x,y,z)\mapsto (x,y,z')$ with $z'=3xy-z$ and we study the quantity $a:=(y-z)(y-z')$; since $z+z'=3xy$ and $zz'=x^2+y^2$, we have 
$$a=2y^2+x^2-3xy^2 = (2y^2-2xy^2)+(x^2-xy^2)<0$$ 
as $2y^2\leq 2xy^2$ and $x^2<xy^2$ thanks to our assumption $1\leq x<y$; because $y-z<0$ (by hypothesis), $a<0$ means that $z'<y<z$. 

In summary, given a Markov triple $(x,y,z)\in\mathbb{N}^3$ with $x\leq y\leq z$, we saw that either $x=y$ and $(x,y,z)\in\{(1,1,1),(1,1,2)\}$, or $x<y<z$ and we can apply a Vieta involution to convert $(x,y,z)$ into a Markov triple $(x,y,z')$ with $x, z'< y < z$. By repeating this argument finitely many times, we see that a sequence of Vieta involutions and permutations of coordinates allows us to convert the Markov triple $(x,y,z)$ into $(1,1,1)$, as desired. 

\subsubsection{The Markov tree} The descent argument above permits us to organize all ordered Markov triples $(x,y,z)\in\mathbb{N}^3$, $x\leq y\leq z$, into the so-called \emph{Markov tree} whose branches connect ordered Markov triples deduced from each other by a Vieta involution (up to permutation of coordinates). 

\begin{figure}[htb!]
\includegraphics[scale=0.4]{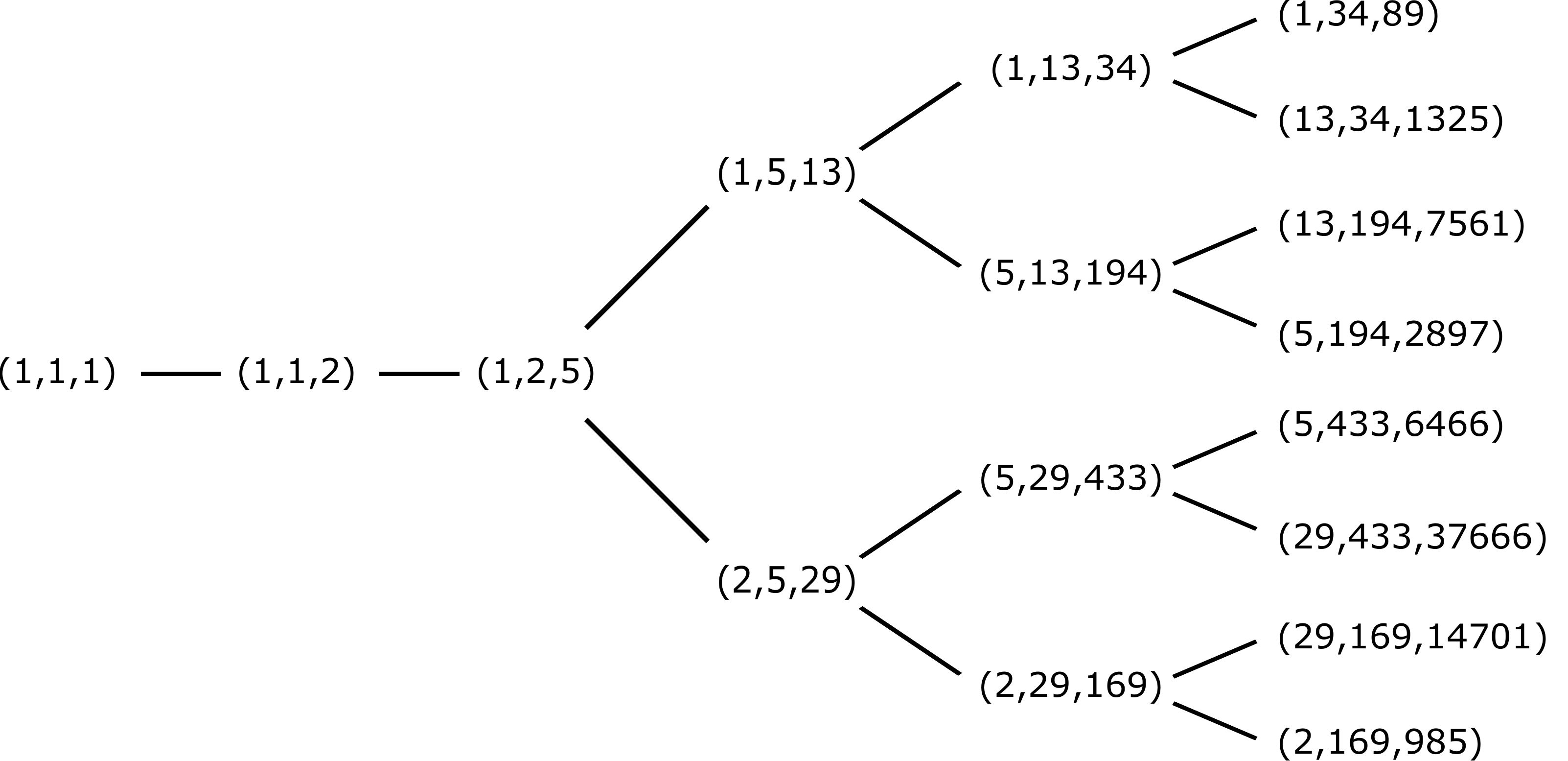}
\end{figure}

The knowledge of Markov's tree permits us to write down the first few elements of $M\cap [\sqrt{5},3)$: since the first few Markov triples are $(1,1,1)$, $(1,1,2)$, $(1,2,5)$, $(1,5,13)$ and $(2,5,29)$, we have that the first few Markov numbers\footnote{The attentive reader certainly noticed that some of these numbers are part of Fibonacci's sequence $(F_n)_{n\in\mathbb{N}}$ and this is not a coincidence: it is possible to check that $(1,F_{2m-1}, F_{2m+1})$ is a Markov triple for all $m\in\mathbb{N}$.} are $1$, $2$, $5$, $13$ and $29$, so that 
$$M\cap[\sqrt{5},3) = \left\{\sqrt{5}, \sqrt{8}, \frac{\sqrt{221}}{5}, \frac{\sqrt{1517}}{13}, \frac{\sqrt{7565}}{29},\dots \right\}$$

\subsubsection{Beyond the Markov tree ...} The Markov tree and numbers are fascinating objects. For instance, it was conjectured by Frobenius in 1913 that Markov triples $(x,y,z)\in\mathbb{N}^3$, $x\leq y\leq z$ are actually determined by the Markov number $z$ (cf. Bombieri's survey article \cite{MR2345177}). 

Also, Zagier \cite{MR669663} showed that the number $M(x)$ of Markov numbers below $x$ is 
$$M(x) = c(\log x)^2+O(\log x (\log\log x)^2)$$ 
where $c=0.180717104711507\dots$ is an \emph{explicit} constant\footnote{Furthermore, Zagier conjectured that $M(x) = c(\log(3x))^2+o(\log x)$, so that $n$-th Markov number $z_n$ is $z_n\sim \frac{1}{3}A^{\sqrt{n}}$ where $A=e^{1/\sqrt{c}}=10.5101504\dots$.}, and, more recently, Baragar \cite{MR1295950} and Gamburd--Magee--Ronan \cite{MR4024562} studied the general problem of counting integral points on the Markov--Hurwitz varieties of the form 
$$x_1^2+\dots+x_n^2 = a x_1\dots x_n +k$$
where $n\geq 3$, $a\geq 1$ and $k$ are integers. 

Moreover, the Markov triples are related to lengths of simple closed geodesics on a certain \emph{hyperbolic} once-punctured torus: in fact, the commutator subgroup $\Gamma$ of $SL(2,\mathbb{Z})$ is an index $12$ subgroup generated by $A_0=\left(\begin{array}{cc}1&-1\\-1&2\end{array}\right)$ and $B_0=\left(\begin{array}{cc}1&1\\1&2\end{array}\right)$; the quotient $SL(2,\mathbb{R})/\Gamma$ is the unit cotangent bundle of a hyperbolic once-punctured torus whose simple closed geodesics correspond to the elements $A\in\Gamma$ in a pair $(A,B)$ of generators of $\Gamma$; the hyperbolic lengths of these geodesics have the form $2\cosh^{-1}(tr(A)/2)$, so that they are related to Markov triples because Fricke proved that any generating pair $(A,B)$ of $\Gamma$ satisfies 
$$tr(A)^2+tr(B)^2+tr(AB)^2 = tr(A) tr(B) tr(AB),$$
i.e., $(tr(A)/3,tr(B)/3,tr(AB)/3)$ is a Markov triple. 

Furthermore, it is known (cf. \cite{MR2026539}) that the level sets of the function $\kappa(x,y,z)=x^2+y^2+z^2-xyz-2$ parametrize the elements of the $SL(2,\mathbb{R})$-character variety\footnote{Naively speaking, the $G$-character variety of a topological surface $S_{g,n}$ of genus $g$ with $n$ punctures is the set of equivalence classes of representations $\rho:\pi_1(S_{g,n})\to G$ modulo the natural action of $G$ by conjugation.} of once-punctured torii and each Markov triple $(x,y,z)$ produces an integral point $(3x,3y,3z)$ of $\kappa^{-1}(-2)$.  

Finally, Bourgain--Gamburd--Sarnak \cite{MR3456887} investigated the family of graphs $(\mathcal{G}_p)$ (indexed by the set of prime numbers $p$) obtained by applying Vieta involutions and permutation of coordinates to the solutions in $\mathbb{F}_p^3\setminus\{(0,0,0)\}$ to Markov--Hurwitz equation $x^2+y^2+z^2=3xyz$. In this setting, they showed $\mathcal{G}_p$ has a giant component $\mathcal{C}_p$ in the sense that $\#(\mathcal{G}_p\setminus\mathcal{C}_p) = O_{\varepsilon}(p^{\varepsilon})$ for all $\varepsilon > 0$, and they used the technology involved in the proof of this statement to establish that almost all Markov numbers are composite, i.e., 
$$\frac{\#\{p \textrm{ prime Markov number}\leq T\}}{\#\{z \textrm{ Markov number}\leq T\}}\to 0$$ as $T\to\infty$. Also, they conjectured that the graphs $\mathcal{G}_p$ are connected\footnote{The connectedness of $\mathcal{G}_p$ for all large $p$ was very recently established by Chen: cf. the arxiv preprint arXiv:2011.12940.} and they form an expander family\footnote{That is, there is an \emph{uniform} spectral gap for the adjacency matrices of these graphs.}. 

\subsection{Continued fractions}

The definition of the Lagrange spectrum suggests that we can study $L$ provided there is a method to find the best rational approximations of a given irrational number $\alpha$ (such as $22/7$ and $355/113$ for $\pi$). 

As it turns out, one can guess the best rational approximations for $\alpha$ out of its \emph{continued fraction expansion}. More precisely, given an irrational number $\alpha_0=\alpha$, let $a_0=\lfloor\alpha\rfloor$, so that $\alpha_0-a_0\in (0,1)$. We define recursively $\alpha_1=\frac{1}{\alpha_0-a_0}$ and $a_n=\lfloor\alpha_n\rfloor\in\mathbb{N}^*$, $\alpha_{n+1}=\frac{1}{\alpha_n-a_n}$ for all $n\in\mathbb{N}^*$. In this context, we say that $\alpha$ has continued fraction expansion 
$$\alpha=a_0+\frac{1}{a_1+\frac{1}{a_2+\frac{1}{\ddots}}} =[a_0; a_1, a_2, \dots]$$ 
and we denote by 
$$\frac{p_n}{q_n} = a_0+\frac{1}{a_1+\frac{1}{\ddots+\frac{1}{a_n}}} =[a_0; a_1, \dots, a_n]$$ 
the \emph{convergents} of $\alpha$. For example, $\pi=[3;7,15,1,292,1,1,1,2,1,3,1,14,2,1,\dots]$, so that $\frac{p_1}{q_1}=22/7$ and $p_3/q_3 = 355/113$. 

It is possible to prove that $(p_n/q_n)$ provides the best rational approximations for $\alpha$ because 
\begin{itemize}
\item if $\frac{p}{q}\in\mathbb{Q}$ and $|\alpha-\frac{p}{q}|<\frac{1}{2q^2}$, then $\frac{p}{q}=\frac{p_n}{q_n}$ for some $n\in\mathbb{N}$; 
\item for each $m\in\mathbb{N}$, $|\alpha-\frac{p_m}{q_m}|<1/q_m^2$ and $\min\{q_m|q_m\alpha-p_m|, q_{m+1}|q_{m+1}\alpha-p_{m+1}|\}<1/2$; 
\item for each $m\in\mathbb{N}$, $|\alpha-p/q|>|\alpha-p_m/q_m|$ for all $p/q\neq p_m/q_m$ with $q<q_{m+1}$.  
\end{itemize}

In particular, the best constant $\ell(\alpha)=\limsup\limits_{p,q\to\infty}\frac{1}{|q(q\alpha-p)|}$ of Diophantine approximation for $\alpha$ depends only on its convergents, i.e., 
$$\ell(\alpha) = \limsup\limits_{n\to\infty}\frac{1}{|q_n(q_n\alpha-p_n)|}.$$

\subsection{Perron's definition of the spectra}\label{ss.Perron} The basic formula 
$$\frac{1}{q_n(q_n\alpha-p_n)} = (-1)^n(\alpha_{n+1}+\beta_{n+1})$$ 
where $\beta_{n+1}:=q_{n-1}/q_n=[0;a_n,\dots, a_1]$ led Perron to propose in 1921 the following \emph{dynamical} interpretation of $L$. Let $\Sigma = (\mathbb{N}^*)^{\mathbb{Z}}$ be the (non-compact) symbolic space of bi-infinite sequences of non-zero natural numbers. The shift map $\sigma:\Sigma\to\Sigma$ is the dynamical system given by 
$$\sigma((a_n)_{n\in\mathbb{Z}}) = (a_{n+1})_{n\in\mathbb{Z}}.$$ 
In this language, the Lagrange spectrum is the set of finite \emph{asymptotic} records of heights of the orbits of $\sigma$ with respect to the (proper) height function $f:\Sigma\to\mathbb{R}$, $f((a_n)_{n\in\mathbb{Z}}) := [a_0; a_1, a_2,\dots] + [0; a_{-1}, a_{-2},\dots]$, i.e., 
$$L=\left\{\limsup\limits_{n\to\infty} f(\sigma^n(x)) < \infty : x\in \Sigma\right\}.$$ 

Interestingly enough, one can use the classical reduction theory of binary quadratic forms (due to Lagrange and Gauss) to prove that the Markov spectrum is the set of finite \emph{absolute} records of heights of the orbits of $\sigma$ with respect to $f$, i.e., 
$$M=\left\{\sup\limits_{n\in\mathbb{Z}} f(\sigma^n(x)) < \infty : x\in \Sigma\right\}.$$

From these dynamical characterizations of $L$ and $M$, Perron deduced that 
\begin{itemize}
\item $\sup\limits_{n\in\mathbb{Z}}f(\sigma^n(x))\leq \sqrt{12}$ if and only if $x\in\{1,2\}^{\mathbb{Z}}$; 
\item $\sqrt{12}, \sqrt{13}, \frac{9\sqrt{3}+65}{22}\in L$;
\item $M\cap (\sqrt{12}, \sqrt{13}) = M\cap(\sqrt{13}, \frac{9\sqrt{3}+65}{22}) = \emptyset$.
\end{itemize}

Moreover, one can use\footnote{In the sequel, \emph{eventually periodic} means eventually periodic on \emph{both} sides (perhaps with different periods).} this dynamical point of view to prove that
$$L=\overline{\{\sup\limits_{n\to\infty} f(\sigma^n(y)): y\in \Sigma \textrm{ is periodic}\}}$$
and 
$$M=\overline{\{\sup\limits_{n\in\mathbb{Z}} f(\sigma^n(z)): z\in \Sigma \textrm{ is eventually periodic}\}}.$$ 
Thus, $L\subset M$ are \emph{closed} subsets of the real line. 

\subsection{Dynamics on the modular surface}\label{ss.modular-surface} The shift map $\sigma:\Sigma\to\Sigma$ can be thought as an invertible map extending the \emph{Gauss map} $G:(0,1]\to [0,1)$, $G(x)=\{1/x\}$. Indeed, the definitions imply that the Gauss map acts on continued fraction expansions by left-shift on \emph{half-infinite} sequences of natural numbers: 
\begin{equation}\label{e.Gauss-shift}
G([0;a_1, a_2,\dots]) = [0; a_2,\dots].
\end{equation}

Using the well-known link (due to Artin, Cohn, Series, Arnoux, ...) between the Gauss map and the geodesic flow $g_t$ on the unit cotangent bundle $SL(2,\mathbb{R})/SL(2,\mathbb{Z})$ to the modular surface (cf. \cite{MR1279059}), one can also describe the Lagrange spectrum as the set of finite asymptotic records of the heights of the orbits of a continuous-time, smooth dynamical system, namely, 
$$L=\{\limsup\limits_{t\to\infty} H(g_t(x))<\infty: x\in SL(2,\mathbb{R})/SL(2,\mathbb{Z})\}$$ 
where $H:SL(2,\mathbb{R})/SL(2,\mathbb{Z})\to\mathbb{R}$ is a certain (proper) function\footnote{By thinking of $SL(2,\mathbb{R})/SL(2,\mathbb{Z})$ as the space of unimodular lattices in $\mathbb{R}^2$, one has $H(x)=2/\textrm{sys}(x)^2$, where $\textrm{sys}(x)$ is the systole of $x\simeq g(\mathbb{Z^2})$, $g\in SL(2,\mathbb{R})$.}.  

\subsection{The end of the classical spectra} The expression of the height function $f:\Sigma\to \mathbb{R}$ in Perron's definition of the spectra suggests that $L$ and $M$ are related to \emph{arithmetic sums} of Cantor sets of real numbers whose continued fraction expansions have restricted digits. 

In other terms, the study of \emph{projections} of products of certain Cantor sets under the function $\pi:\mathbb{R}^2\to \mathbb{R}$, $\pi(x,y):=x+y$, should provide some insights into the fine structures of $L$ and $M$.  

This idea was explored by Hall in 1947 to show that $L$ contains the half-line $[6,\infty)$. For this sake, Hall considered the Cantor set $C(4)=\{[0;a_1,a_2,\dots]:1\leq a_i\leq 4 \,\,\, \forall i\in\mathbb{N}\}$ and he established that 
\begin{eqnarray*}
C(4)+C(4)&:=&\{x+y: (x,y)\in C(4)\times C(4)\} \\ 
&=& [\sqrt{2}-1, 4(\sqrt{2}-1)]
\end{eqnarray*} 
is an interval of length $>1$. This fact implies that given $\ell\in[6,\infty)$, one can find $c_0\in\mathbb{N}$ such that $5\leq c_0\leq \ell$ and $\ell-c_0\in C(4)+C(4)$, say 
$$\ell=c_0+[0;a_1,a_2,\dots]+[0;b_1,b_2,\dots]$$ 
with $1\leq a_i, b_i\leq 4$ for all $i\in\mathbb{N}$. Thus, the irrational number $\alpha$ with continued fraction expansion 
$$\alpha=[0;\underbrace{b_1,c_0,a_1}_{1^{st } block}, \dots, \underbrace{b_n,\dots, b_1, c_0, a_1, \dots, a_n}_{n^{th } block}, \dots]$$ 
satisfies $\ell=\ell(\alpha)\in L$. Since $\ell\geq 6$ was arbitrary, we conclude that $L\supset [6,\infty)$. 

The largest half-line of the form $[c_F,\infty)$ included in $L$ is called \emph{Hall ray}. In 1975, Freiman famously claimed that 
$$c_F = 4+\frac{253589820 + 283798\sqrt{462}}{491993569} = 4.527829566\dots$$ 

\subsection{Intermediate portions of $L$ and $M$}\label{ss.inter} We saw above that $L\subset M$ are closed subsets of the real line such that $L\cap (-\infty, 3)=M\cap (-\infty, 3) = \{\sqrt{5}, \sqrt{8}, \dots\}$, $L\cap[c_F,\infty)=M\cap[c_F,\infty) = [c_F,\infty)$, and $\sqrt{12}, \sqrt{13}\in L$, but $(\sqrt{12},\sqrt{13})\cap M=\emptyset$. 

In particular, $L$ and $M$ \emph{coincide} in the portions $(-\infty, 3)$, $[\sqrt{12}, \sqrt{13}]$, and $[c_F,\infty)$. Nonetheless, it was discovered by Freiman in 1968 that $M\setminus L\neq\emptyset$: more concretely, Freiman found a \emph{countable} subset of $M\setminus L$ located near $3.11$. 

Subsequently, Freiman discovered (in 1973) a new element of $M\setminus L$ near $3.29$, and Flahive showed\footnote{Actually, in his unpublished PhD thesis from 1976, Y.-C. You found an \emph{uncountable} subset of $M\setminus L$ near $3.29$ which is bi-Lipschitz homeomorphic to the Cantor set of continued fraction expansions obtained by concatenations of the words $11$ and $22$.} in 1977 that this element is the limit of an explicit sequence of elements of $M\setminus L$ near $3.29$. 

This ``concentration'' of examples of elements of $M\setminus L$ between $3$ and $3.3$ led Cusick to \emph{conjecture} in 1975 that $L$ and $M$ should coincide above $\sqrt{12}$. 

In any case, the previous paragraphs hint that the intermediate portions of the spectra (between $3$ and $c_F$) might have a \emph{complicated} structure in comparison to their beginning and ending. 

In 1971, Hall gave an upper bound on the \emph{fractal complexity} of $M\cap[\sqrt{5},\sqrt{10}]$. More precisely, he used Perron's definition of the spectra to establish that $M\cap[3,\sqrt{10}]\subset 2+ U + U$ where 
$$U=\{[0;a_1,a_2,\dots]: (a_i a_{i+1} a_{i+2})\neq (121) \, \, \forall \, i\in\mathbb{N}^*\}.$$ 
After that, he analysed the sizes of all intervals covering $U$ with extremities of the form $[0;a_1,\dots, a_n, \overline{2122}]$ and $[0;a_1,\dots, a_n,\overline{1222}]$ (where $\overline{w}$ stands for the periodic sequence obtained by infinite concatenation of a string $w$) in order to show that the \emph{Hausdorff dimension}\footnote{Recall that the Hausdorff dimension $0\leq \textrm{dim}(X)\leq m$ of $X\subset\mathbb{R}^m$ measures how hard it is to efficiently cover $X$: in fact, $\textrm{dim}(X)=\inf\left\{s>0: \inf\limits_{X\subset \bigcup\limits_{n\in\mathbb{N}} B(x_n, r_n)} \sum\limits_{n\in\mathbb{N}} r_n^s =0\right\}$. Thus, a countable set has zero Hausdorff dimension and any set $X\subset \mathbb{R}^m$ with $\textrm{dim}(X)<m$ has zero $m$-dimensional Lebesgue measure.} of $U$ is at most $0.465$. From this estimate, it is possible to infer that the Hausdorff dimension of $U+U = \pi(U\times U)$ has Hausdorff dimension at most $2\cdot\textrm{dim}(U) < 0.93$ and, \emph{a fortiori}, $M\cap[\sqrt{5},\sqrt{10}]\subset 2 + U + U$ has zero Lebesgue measure. 

On the other hand, it is believed that $L$ and $M$ should contain non-trivial intervals between $t_1$ and $c_F$: for instance, a folklore question (appearing in page 71 of Cusick--Flahive book \cite{MR1010419}) asks whether $L\cap[\sqrt{5},\sqrt{12}]$ has non-empty interior, and Berstein conjectured in 1973 that $[4.1, 4.52]\subset L$. Here, it is worth pointing out that the inspiration for the first (folklore) question comes from: 
\begin{itemize}
\item Perron's result that $M\cap[\sqrt{5}, \sqrt{12}]$ is closely related to the arithmetic sum $C(2)+C(2)$ where $C(2):=\{[0;a_1,a_2,\dots]: 1\leq a_i\leq 2\,\,\forall\,i\}$; 
\item the expectation\footnote{The typical projections of generic planar Cantor sets with Hausdorff dimension $>1$ are expected to contain intervals thanks to the combination of Marstrand's theorem asserting that the projections in almost all directions of a subset of $\mathbb{R}^2$ with Hausdorff dimension $>1$ have positive Lebesgue measures and the work of Yoccoz and the second author on the stable intersections of Cantor sets via renormalization techniques.} that $C(2)+C(2)$ contains intervals because $C(2)+C(2)$ is the projection $\pi(C(2)\times C(2))$ of a planar ``non-linear'' Cantor set $C(2)\times C(2)\subset\mathbb{R}^2$ with Hausdorff dimension $2\cdot\textrm{dim}(C(2)) > 1$. 
\end{itemize}
Also, Berstein thinks that $[4.1, 4.52]\subset L$ because of Freiman's work on the computation of the beginning $c_F$ of Hall's ray. 

\subsection{Recent results about $M\cap(3,c_F)$}\label{ss.main-results} Despite the strong belief (expressed by the conjectures and questions in the previous subsection) that $M\cap (3, c_F)$ must have an intricate structure, the first rigorous result in this direction was obtained only in 2018 by the second author \cite{MR3815461}. In fact, he showed that: 
\begin{itemize}
\item for each $t\in\mathbb{R}$, the Hausdorff dimension $d(t)$ of $L\cap (-\infty, t)$ coincides with the Hausdorff dimension of $M\cap (-\infty, t)$, i.e., $M\setminus L$ is \emph{not} big enough to create jumps in dimension between $L$ and $M$; 
\item $d(t)$ is a \emph{continuous}, non-H\"older function of $t$ such that $d(3+\varepsilon) > 0$ for all $\varepsilon>0$ and $d(\sqrt{12})=1$.   
\end{itemize} 
The second item above implies that $L$ and $M$ necessarily contain complicated fractal sets. For instance, it is not difficult to check that if $K\subset\mathbb{R}$ is a Cantor set defined by simple interactive (dynamical) rules like Cantor's middle-third set\footnote{Recall that Cantor's middle-third set is given by breaking $[0,1]$ into three intervals of equal lengths, removing the subinterval in the middle, and then successively repeating this operation to each of the remaining subintervals.}, then the function $t\mapsto \textrm{dim}(K\cap (-\infty, t))$ is always a piecewise constant and discontinuous: 
$$\textrm{dim}(K\cap (-\infty, t)) = \left\{ \begin{array}{cc} 0 & \textrm{ if } t \leq \min K \\ \textrm{dim}(K) & \textrm{ if } t > \min K\end{array}\right.$$

More recently, we investigated in \cite{MR3860483}, \cite{MR3925086} and \cite{MR4152626} the fine structure of $M\setminus L$ and we proved that it is \emph{richer} than conjectured by Cusick. More precisely, there are three open intervals $I_1$, $I_2$ and $I_3$ near $3.11$, $3.29$ and $3.7$ such that: 
\begin{itemize}
\item the sizes of $I_1, I_2, I_3$ are $\sim 2\cdot 10^{-10}, 2\cdot 10^{-7}, 10^{-10}$; 
\item the extremities of $I_j$ belong to $L$, but $L\cap I_j=\emptyset$ for each $1\leq j\leq 3$; 
\item $(M\setminus L)\cap I_j$, $1\leq j\leq 3$, are \emph{closed} subsets with Hausdorff dimensions $>0.26, 0.353, 0.531$ resp.; 
\item $I_1$ and $I_2$ contain the examples of elements of $M\setminus L$ previously found by Freiman and Flahive, and the elements of $(M\setminus L)\cap I_3$ provide a \emph{negative} answer to Cusick's conjecture that $L$ and $M$ should coincide above $\sqrt{12}$. 
\end{itemize} 
On the other hand, we proved that $M\setminus L$ is \emph{not} very rich\footnote{We also offered some \emph{heuristic} evidence towards an upper bound of the form $\textrm{dim}(M\setminus L)<0.888$ and, as it turns out, Pollicott--Vytona proved that $\textrm{dim}(M\setminus L)<0.8822195$ in a preprint available at arXiv:2012.07083.} because $\textrm{dim}(M\setminus L)<0.987$.\footnote{Very recently, in collaboration with M. Pollicott and P. Vytnova, in a preprint available at arXiv:2106.06572, we improved these estimates, showing that $0.537152 < \textrm{dim}(M\setminus L) < 0.796445$.} 

\begin{figure}[htb!]
\includegraphics[scale=1.4]{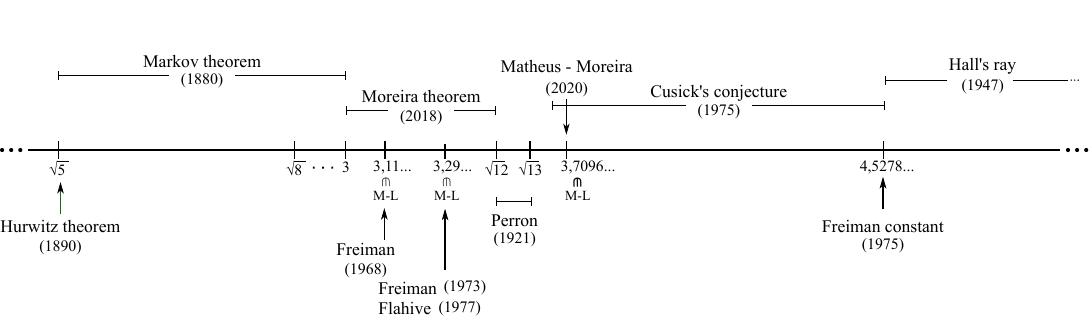}
\end{figure}

Besides the metric results discussed above, our techniques also give \emph{topological} consequences for $L$ and $M$: for example, the subset $L'$ of non-isolated points of $L$ is \emph{perfect}, i.e., $L'=L''$, and the \emph{interiors} of $L$ and $M$ coincide. However, we \emph{ignore} whether $M'$ is a perfect set, and, contrary to the initial impression given by the fact that $(M\setminus L)\cap I_j$, $1\leq j\leq 3$, are closed subsets, Lima, Vieira and the authors proved that $M\setminus L$ is \emph{not} closed\footnote{We also gave some \emph{evidence} towards the possibility that $3\in L\cap\overline{(M\setminus L)}$.} by establishing that $1+3/\sqrt{2}\in L\cap\overline{(M\setminus L)}$. 

In an attempt to further investigate interesting questions about the structure of $M\cap (3, c_F)$, Delecroix and the authors \cite{MR4109576} developed an algorithm providing $1/Q$-approximations\footnote{Here, we mean close in \emph{Hausdorff topology}, i.e., $A, B\subset\mathbb{R}$ are $\delta$-close if and only if for each $a\in A$ and $b\in B$ there are $c\in B$ and $d\in A$ such that $|a-c|, |b-d|\leq \delta$.} to $L$ and $M$ after a running time $O(Q^{2.367})$. This algorithm was implemented (on Sage) by Delecroix to produce the figure below of $L_2:=L\cap [\sqrt{5}, \sqrt{12}]$, but unfortunately, we could not use this algorithm yet to get definite ideas about Berstein's conjecture that $[4.1, 4.52]\subset L$. Nevertheless, we hope that some variant of this algorithm will be helpful in the future because its running time is not very big in comparison\footnote{In fact, a back-of-the-envelope calculation reveals that the naive algorithm requires computations with continued fractions associated to $\sim 4^{Q^4}$ strings of lengths $\leq Q^4$ to rigorously produce a $1/Q$-approximation to $L$ and $M$.} with the ``naive'' algorithm stemming from the characterisations of $L$ and $M$ via the closures of the values of height records of periodic and eventually periodic elements of $\Sigma=(\mathbb{N}^*)^{\mathbb{Z}}$. 

\begin{figure}[htb!]
\includegraphics[scale=0.7]{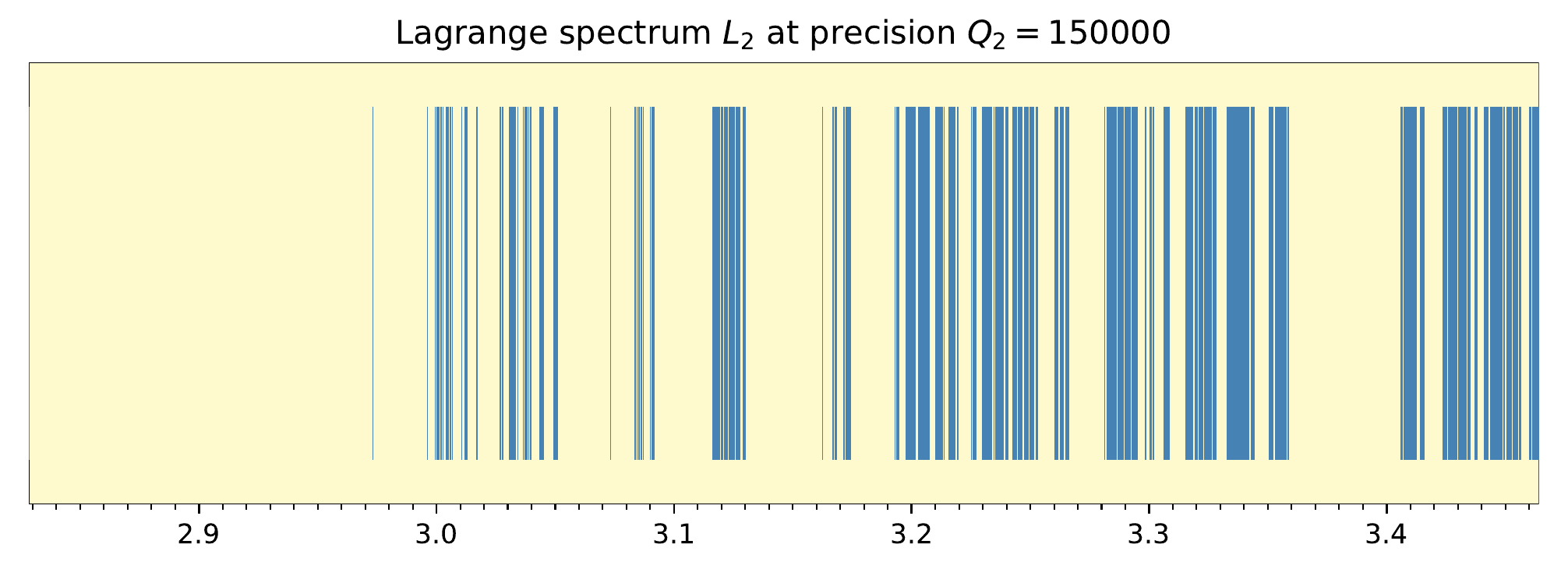}
\end{figure}

The common theme behind all results described in this subsection is the study of portions of $L$ and $M$ via the fractal geometry of certain types of \emph{dynamically defined Cantor sets}. In the next section, we will briefly discuss these objects and we will make some comments about their applications to the proofs of our results on the Hausdorff dimensions of $L\cap(-\infty, t)$, $M\cap (-\infty, t)$ and $M\setminus L$. 

\section{Dynamical Cantor sets} A \emph{dynamically defined} Cantor set $K\subset R$ is 
$$K=\bigcap\limits_{n\in\mathbb{N}} \psi^{-n}(I_1\cup\dots\cup I_k),$$ 
where $\psi: I_1\cup\dots\cup I_k\to I$ is a smooth map from a collection of disjoint compact intervals $I_1,\dots, I_k$ onto the convex hull $I$ of their union such that 
\begin{itemize}
\item $\psi$ is \emph{expanding}: $|\psi'(x)|>1$ $\forall$ $x\in I_1\cup\dots\cup I_k$; 
\item $\{I_1, \dots, I_k\}$ is a \emph{Markov partition}: each $\psi(I_i)$ is the convex hull of the union of some of the intervals $I_j$; 
\item $\psi$ is \emph{topologically mixing}: there exists $n_0\in\mathbb{N}$ with $\psi^{n_0}(K\cap I_i)=K$ for all $1\leq i\leq k$. 
\end{itemize} 

The most famous example of dynamically defined Cantor set is arguably Cantor's middle-third set 
$$K_{1/3} = \left\{\sum\limits_{n=1}^{\infty}\frac{a_n}{3^n}: (a_n)_{n\in\mathbb{N}}\in \{0,2\}^{\mathbb{N}}\right\}$$ 
obtained by dividing $[0,1]$ into three intervals with equal lengths, removing the middle one $(1/3, 2/3)$, dividing the remaining intervals into three equal subintervals, removing the middle ones, repeating this procedure \emph{ad infinitum}, and keeping only the points of $[0,1]$ never falling into an excluded subinterval. 

In fact, $K_{1/3}=\bigcap\limits_{n\in\mathbb{N}} \phi^{-n}(I_1\cup I_2)$ where $I_1=[0,1/3]$, $I_2=[2/3,1]$, and 
$$\phi(x)=\left\{ \begin{array}{cc} 3x & \textrm{ if } x\in I_1, \\ 3x-2 & \textrm{ if } x\in I_2. \end{array}\right.$$
  
The geometry of dynamical Cantor sets was studied in depth by several authors over the last 100 years. In particular, we dispose nowadays of several methods to compute the Hausdorff dimension of dynamical Cantor sets. For instance, given a dynamical Cantor set $K$ generated by $\psi:I_1\cup\dots\cup I_k\to I$, for each $m\in\mathbb{N}$, we can consider the collection $\mathcal{R}_m$ of intervals consisting of the connected components of $\bigcap\limits_{n=1}^m\psi^{-n}(I_1\cup\dots\cup I_k)$. By definition, $\mathcal{R}_m$ is a cover of $K$ and it can be used to compute $\textrm{dim}(K)$: if we denote by $\lambda(J)=\min|\psi'|_J|$, $\Lambda(J)=\max|\psi'|_J|$ and we fix a mixing time $n_0$ for $\psi$ (i.e., $\psi^{n_0}(I_j\cap K)=K$ for each $1\leq j\leq k$), then it is possible to prove that $\alpha_m\leq \textrm{dim}(K)\leq \beta_m$, where 
$$\sum\limits_{J\in\mathcal{R}_m} \frac{1}{\lambda(J)^{\beta_m}} = 1 \textrm{ and} \sum\limits_{J\in\mathcal{R}_m} \frac{1}{\Lambda(J)^{\alpha_m}} = \max|(\psi^{n_0-1})'|$$ 
(cf. pages 68 to 70 of Palis--Takens book \cite{MR1237641}). These bounds allow an \emph{exact} calculation of the Hausdorff dimension of dynamical Cantor sets associated to piecewise \emph{affine} maps $\psi$ with full branches (i.e., mixing time $n_0=1$) such as Cantor's middle-third set $K_{1/3}$: indeed, $K_{1/3}$ is defined by a map $\phi:I_1\cup I_2\to I$ such that $\phi'\equiv 3$ (and $n_0=1$), so that $\alpha_1\leq \textrm{dim}(K_{1/3})\leq \beta_1$, where 
$$2(1/3)^{\alpha_1} = 1 = 2(1/3)^{\beta_1},$$ 
i.e., $\textrm{dim}(K_{1/3})=\log 2/\log 3$. 

Unfortunately, the elementary technique described in the previous paragraph doesn't permit us to compute the dimension of dynamical Cantor sets $K$ associated to \emph{non-essentially affine}\footnote{In the sense that there is no smooth change of coordinates $h:I\to J$ making $h\circ \psi\circ h^{-1}$ into a piecewise affine map.} maps $\psi:I_1\cup\dots\cup I_k\to I$. In fact, one can check that $\beta_m-\alpha_m=O(1/m)$ in general, so that the convergence of $\alpha_m$ and $\beta_m$ to $\textrm{dim}(K)$ might be quite slow\footnote{Note that the calculations of $\alpha_m$ and $\beta_m$ require us to manipulate the intervals in $\mathcal{R}_m$ and, hence, our $O(1/m)$-approximation to $\textrm{dim}(K)$ comes from the computations with a quantity $\#\mathcal{R}_m$ of intervals growing exponentially with $m$.} in the non-essentially affine situations related to the classical spectra $L$ and $M$.

Nevertheless, Bowen \cite{MR556580} discovered in 1979 a famous formula for the Hausdorff dimension of dynamical Cantor sets $K$ which was subsequently explored by several authors (including Falk, Hensley, Jenkinson, McMullen, Nussbaum, Pollicott, Vytnova) for a \emph{fast} computation of \emph{several} digits of $\textrm{dim}(K)$. Roughly speaking, Bowen's formula for a dynamical Cantor set $K$ associated to an expanding map $\psi:I_1\cup\dots\cup I_k\to I$ starts with a family $(\mathcal{L}_t)_{t\in (0,1)}$ of \emph{Ruelle--Perron--Frobenius transfer operators} 
$$\mathcal{L}_t f (x) = \sum\limits_{y\in\psi^{-1}(x)} f(y) |\psi'(y)|^{-t}$$ 
acting on adequate spaces of smooth functions $f:I\to \mathbb{R}$. This kind of operator was originally introduced by Ruelle in his study of ergodic theoretical properties of dynamical systems via an analogy with the so-called \emph{thermodynamical formalism} in statistical mechanics.\footnote{Transfer operators describe the action of $\psi$ on the densities (Radon--Nykodym derivatives) of probability measures: by fixing a probability measure $\lambda$ which is non-singular for $\psi$ (in the sense that $\psi_*(\lambda) = \lambda\circ \psi^{-1}$ is absolutely continuous with respect to $\lambda$), we get a transfer operator $\mathcal{L}_{(\lambda)}$ by looking at the Radon--Nykodym derivative $\mathcal{L}_{(\lambda)}(f) = \frac{d\psi_*(f\lambda)}{d\lambda}$ of the push-forward $\psi_*(f\lambda)$ under $\psi$ of a probability $\nu = f \lambda$ which is absolutely continuous with respect to $\lambda$.} The operators $\mathcal{L}_t$ are \emph{quasi-compact} (i.e., their spectral theories share many parallels with matrices in finite-dimensional vector spaces) and they possess a leading eigenvalue $\lambda_t$ with multiplicity one. In this setting, Bowen showed that the Hausdorff dimension of $K$ is the unique parameter $\textrm{dim}(K)\in (0,1)$ such that 
$$\lambda_{\textrm{dim}(K)}=1.$$ 
In particular, Cantor's middle-third set $K_{1/3}$ has a Hausdorff dimension $\log 2/\log 3$ such that $\mathcal{L}_{\log 2/\log 3}(1)=1$. 

The spectral characterisation of $\textrm{dim}(K)$ provided by Bowen can be reformulated by saying that $s=\textrm{dim}(K)$ is the unique parameter such that 
$$\det(\textrm{Id}-\mathcal{L}_s)=0.$$ 
In other words, we are interested in the value of $s$ such that $z=1$ is a zero of $\det(\textrm{Id}-z\mathcal{L}_s)$. As it turns out, we can efficiently compute this value thanks to the \emph{Fredholm determinant expansion}\footnote{This expansion generalizes the case of finite-dimensional vector space: if $K$ is a $n\times n$ matrix with eigenvalues $\lambda_1,\dots, \lambda_n$, then $\det(\textrm{Id}+K) = \prod\limits_{j=1}^n(1+\lambda_i)$, so that $\log\det(\textrm{Id}+K) = \sum\limits_{j=1}^n\log(1+\lambda_i) = \sum\limits_{m=1}^{\infty}\frac{(-1)^{m+1}}{m}\left(\sum\limits_{j=1}^n\lambda_j^m\right) = \sum\limits_{m=1}^{\infty}\frac{(-1)^{m+1}}{m}\textrm{tr}(K^m)$.} 
$$\det(\textrm{Id}-z\mathcal{L}_s) = \exp\left(-\sum\limits_{n=1}^{\infty}\textrm{tr}(\mathcal{L}_s^n) \frac{z^n}{n}\right).$$
Indeed, it is possible to check that the traces $\textrm{tr}(\mathcal{L}_s^n)$ have nice expressions in terms of the \emph{periodic points} of $\psi$: 
$$\textrm{tr}(\mathcal{L}_s^n) = \sum\limits_{\psi^n(p)=p} \frac{|(\psi^n)'(p)|^{-s}}{1-(\psi^n)'(p)^{-1}}.$$ 
Thus, we can write $\det(\textrm{Id}-z\mathcal{L}_s) = 1+\sum\limits_{n=1}^{\infty}d_n(s) z^n$, where $d_n(s)$ are \emph{explicit} functions of $s$ and the periodic points of $\psi$. This provides a practical scheme to compute $\textrm{dim}(K)$ when $\psi$ is piecewise real-analytic because it is possible to prove that in this situation the solutions $s_M$ of 
the truncations $1+\sum\limits_{n=1}^{M}d_n(s) = 0$ of Bowen's formula  
$$\det(\textrm{Id}-\mathcal{L}_{\textrm{dim}(K)})=1+\sum\limits_{n=1}^{\infty}d_n(\textrm{dim}(K)) = 0$$
converge quickly to $\textrm{dim}(K)$, namely, 
$$|s_M-\textrm{dim}(K)|=O(\theta^{M^2})$$ 
for some $\theta<1$. 

In particular, the method in the previous paragraph was successfully explored by Jenkinson and Pollicott in 2018 to compute the first 100 decimal digits of the Hausdorff dimension $\textrm{dim}(C(2))$ of the Cantor set $C(2)=\{[0;a_1,a_2,\dots]: 1\leq a_i\leq 2\,\,\forall\,i\}$ which is dynamically defined by the piecewise real-analytic map given by the restriction of the Gauss map to the intervals $[[0;2\overline{12}], [0;2,\overline{21}]$ and $[[0;1\overline{12}], [0;1\overline{21}]$. The outcome of their calculations is that 
$$\textrm{dim}(C(2)) = 0.531280506277205141624468647368\dots$$

\subsection{Dimension of Gauss--Cantor sets} For our purpose of studying the classical spectra, the relevant class of dynamical Cantor sets are the so-called \emph{Gauss--Cantor sets} defined as follows. 

Let $B\subset \bigcup\limits_{n\geq 1}(\mathbb{N}^*)^n$ be a finite set of finite words which is {\it primitive} in the sense that none of its elements is a prefix of another one. The corresponding \emph{complete Gauss--Cantor set} is
$$K(B)=\{[0;\beta_1,\beta_2,\dots]:\beta_i\in B \,\forall i\geq 1\}.$$

The simplest examples\footnote{Note that $C(4)$ and $C(2)$ already appeared in our discussions of the ending and the intermediate portions of the classical spectra $L$ and $M$.} of complete Gauss-Cantor sets are the sets 
$$C(k):=\{[0;a_1,a_2,\dots]: 1\leq a_i\leq k\,\,\forall\,i\}$$ 
for $k\ge 2$.

In general, a \emph{Gauss-Cantor set} is a set of the type
$$K(\gamma,B)=\{[0;\gamma,\beta_1,\beta_2,\dots]:\beta_i\in B \,\forall i\geq 1\},$$
where $\gamma\in\bigcup\limits_{n\geq 1}(\mathbb{N}^*)^n$ is a finite word.

Notice that $K(\gamma,B)$ is the image of $K(B)$ under the bi-Lipschitz homeomorphism 
$$[0;\gamma, x]\mapsto [0;x] = G^{|\gamma|}([0;\gamma,x]),$$ where $G$ is the Gauss map and $|\gamma|$ is the size of $\gamma$. In particular, $K(B)$ and $K(\gamma,B)$ have the same Hausdorff dimension. 

Complete Gauss-Cantor sets are dynamical Cantor sets defined by iterates of the Gauss map $G$. Indeed, $K(B)$ is a dynamical Cantor set $\psi:\bigcup\limits_{\beta\in B} I(\beta)\to I $ where $\psi|_{I(\beta)}=G^{|\beta|}$ and $I(\beta)$ are intervals with extremities of the form $[0;\beta, x_{\beta}]$ and $[0;\beta, y_{\beta}]$ for adequate choices of $x_{\beta}, y_{\beta}\in B^{\mathbb{N}}$. 

Similarly, one can verify that, in general, Gauss-Cantor sets are also dynamical Cantor sets. 

It is not difficult to show that any Gauss--Cantor set is \emph{non-essentially affine}, i.e., it is dynamically defined by a map $\psi:I_1\cup\dots\cup I_k\to I$ such that there is no smooth change of coordinates $h$ making the second derivative of $h\circ \psi \circ h^{-1}$ to vanish identically on $h(K)$. Thus, a Gauss--Cantor set is geometrically more intricate than dynamical Cantor sets given by piecewise affine maps (such as Cantor's middle-third set) and, hence, we do not expect to get an \emph{exact formula} for its Hausdorff dimension (but only some high precision approximation coming from variants of Bowen's formula, for example). 

Nonetheless, the second author discovered\footnote{A version of this formula was obtained by Hochman and Shmerkin.} that the renormalization techniques (including the so-called \emph{scale recurrence lemma}) introduced by Yoccoz and him \cite{MR1865980} in their study of stable intersections of dynamical Cantor sets can be used to prove that if $K$ is a non-essentially affine dynamical Cantor set and $K'$ is an arbitrary dynamical Cantor set, then the projection $\pi(K\times K')=K+K'$ of $K\times K'$ under $\pi(x,y)=x+y$ has \emph{expected} Hausdorff dimension 
$$\textrm{dim}(K+K') = \min\{1, \textrm{dim}(K)+\textrm{dim}(K')\}.$$ 
In particular, Gauss--Cantor sets satisfy the following \emph{dimension formula}: for any finite words $\gamma, \gamma'$ and finite sets of finite words $B, B'$, one has 
$$\textrm{dim}(K(\gamma,B)+K(\gamma', B'))=\min\{1, \textrm{dim}(K(B))+\textrm{dim}(K(B'))\}.$$

Also, it is worth pointing out the following\footnote{The dynamical explanation for this symmetry is the fact that the Gauss map has a smooth, area-preserving, natural extension.} useful ``symmetry'' on the Hausdorff dimension of Gauss--Cantor sets. Given a finite set of finite words $B$, denote by $B^T=\{\beta^T:\beta\in B\}$ the \emph{transpose} of $B$, where $\beta^T:=(a_n,\dots, a_1)$ stands as usual for the transpose of $\beta=(a_1,\dots,a_n)$. Then, the Hausdorff dimension of the Gauss--Cantor sets associated to $B$ and $B^T$ are equal: 
$$\textrm{dim}(K(B)) = \textrm{dim}(K(B^T)).$$ 
In a certain sense, the proof of this fact goes back to Euler: indeed, he proved that for any finite word $\beta$, if $[0;\beta]=p_n/q_n$, then $[0;\beta^T]=r_n/q_n$. Since the lengths of the intervals $I(\beta_1 \beta_2 \dots \beta_k)$ in the $k$-th step of the construction of $K(B)$ depend only on the denominators of the convergents of $[0;\beta_1 \beta_2 \dots \beta_k]$, Euler's result says that $K(B)$ and $K(B^T)$ are Cantor sets constructed from small intervals with comparable lengths, and, \emph{a fortiori}, they have the same Hausdorff dimension. 

\subsection{Dimension across the spectra} As we said in the end of \S\ref{ss.main-results}, we want to study portions of $L$ and $M$ via dynamical Cantor sets. In particular, a central idea towards the main theorem of \cite{MR3815461} about the continuity of the Hausdorff dimension across the classical spectra is to approximate $L\cap (-\infty,t)$ and $M\cap (-\infty,t)$ from inside and outside by arithmetic sums of Gauss--Cantor sets. 

More precisely, let $D(t)=\textrm{dim}(\ell^{-1}(-\infty,t))$ (where $\ell(\alpha)$ is the best constant\footnote{Cf. the beginning of \S\ref{s.classical-spectra}.} of Diophantine approximation of $\alpha$). Given $t$ such that $D(t)>0$ and $\epsilon>0$, the second author proved the existence of 
\begin{itemize} 
\item a parameter $\delta>0$, 
\item a finite set of positive integers $(a_j)_{1\leq j\leq m+1}$, 
\item a finite set of finite prefixes $\{\gamma_j\}_{j=1}^{2m+2}\subset \bigcup\limits_{n\geq 1}(\mathbb{N}^*)^n$, and 
\item two finite sets $B, B'$ of finite words 
\end{itemize} 
such that the translated arithmetic sum of Gauss--Cantor sets
$$a_{m+1}+K(\gamma_{2m+1},B')+K(\gamma_{2m+2},(B')^T)$$ 
is contained in $L\cap (-\infty,t-\delta)$, the union of translated arithmetic sums of Gauss--Cantor sets 
$$\bigcup_{j\leq m}(a_j+K(\gamma_{2j-1},B)+K(\gamma_{2j},B^T))$$ 
contains $M\cap (-\infty,t+\delta)$, and 
$$D(t)-\epsilon<\textrm{dim}(K(B'))\le D(t) \leq \textrm{dim}(K(B))<D(t)+\epsilon.$$ 
In particular, these facts together with the dimension formula and Euler's symmetry imply that 
\begin{eqnarray*}
\min\{1, 2D(t)-2\epsilon\} &\leq& \min\{1, 2 \textrm{ dim}(K(B'))\} \\ &=& \min\{1, \textrm{dim}(K(B'))+\textrm{dim}(K((B')^T)\} \\ &\leq& \textrm{dim}(L\cap (-\infty,t-\delta)) \\ 
&\leq& \textrm{dim}(L\cap (-\infty,t))
\end{eqnarray*} 
and 
\begin{eqnarray*} 
\textrm{dim}(M\cap (-\infty,t)) &\leq& \textrm{dim}(M\cap (-\infty,t+\delta)) \\ &\leq& \min\{1, \textrm{dim}(K(B))+\textrm{dim}(K(B^T)\} \\ &=&\min\{1, 2\textrm{ dim}(K(B))\} \\ &\leq& \min\{1, 2D(t)+2\epsilon)\}.
\end{eqnarray*} 
Since $\epsilon>0$ was arbitrary, we conclude that 
$$\textrm{dim}(L\cap (-\infty,t))=\textrm{dim}(M\cap (-\infty,t))=\min\{1, 2D(t)\}$$ is a continuous function of $t$. 

In 1982, Bumby extended Hall's ideas in \S\ref{ss.inter} to give a computer-assisted argument \emph{indicating} that 
$D(3.33437) < 1/2$ and  $D(3.3344)>1/2$, so that 
$$3.33437<\inf\{t: \textrm{dim}(M\cap[\sqrt{5},t])=1\}<3.3344$$ in view of our current discussion.\footnote{Very recently, in collaboration with M. Pollicott and P. Vytnova - see arXiv:2106.06572, we \emph{rigorously} improved upon Bumby's claim, showing that $t_1=3.334384\dots$, where this value is rigorously accurate to the $6$ decimal places presented.} 

\subsection{Bounds on $\textrm{dim}(M\setminus L)$} 
In \cite{MR4152626}, we explored the geometry of the intersections of the so-called local \emph{stable} and \emph{unstable} sets\footnote{The local stable set of the left-shift map $\sigma:(\mathbb{N}^*)^{\mathbb{Z}}\to (\mathbb{N}^*)^{\mathbb{Z}}$ at $x=(x_n)_{n\in\mathbb{Z}}$ is the subset of $y=(y_n)_{n\in\mathbb{Z}}$ such that $y_m=x_m$ for all $m\geq 0$. This nomenclature is justified by the fact that the distance between $\sigma^k(y)$ and $\sigma^k(x)$ goes exponentially fast to zero as $k\to+\infty$ when $y$ belongs to the local stable set of $x$. Anagously, we get local unstable sets after replacing $\sigma$ by $\sigma^{-1}$ in the previous definition.} of the shift dynamics $\sigma$ to prove that $M\setminus L \supset a+K(\gamma,B)$, where $a=[3;3,\overline{2,1,2,2,2,3,3}]$, $\gamma=(2,2,2,1,2,3,3,2,2,2,1,2,2,1,2,1,2,1,2)$ and $B=\{1,2\}$. Thus, $\textrm{dim}(M\setminus L)\geq \textrm{dim}(K(\gamma,B))=\textrm{dim}(K(B))=\textrm{dim}(C(2))>0.531.$

After this brief discussion of lower bounds on $\textrm{dim}(M\setminus L)$, let us now sketch\footnote{These ideas can be adapted to give nontrivial upper estimates for $\textrm{dim}(M\setminus L)$ by analyzing other regions of the spectra.} a proof of the upper bound $\textrm{dim}((M\setminus L)\cap [\sqrt{5}, \sqrt{12}]) < 0.93$. 
 
We saw in \S\ref{ss.Perron} that if $\theta\in(\mathbb{N}^*)^{\mathbb{Z}}$, then 
$$m(\theta):=\sup\limits_{n\in\mathbb{Z}}f(\sigma^n(\theta))\leq\sqrt{12} \iff \theta\in\{1,2\}^{\mathbb{Z}}.$$ 
Moreover, we saw in \S\ref{ss.inter} that $m(x)\leq\sqrt{10}$ implies that $x\in\{1,2\}^{\mathbb{Z}}$ doesn't contain $121$ and, \emph{a fortiori},  
$$\textrm{dim}(M\cap[3,\sqrt{10}]) < 0.93.$$ 

Since $(M\setminus L)\cap [3,\sqrt{12}]\subset \left(M\cap[3,\sqrt{10}]\right) \cup \left((M\setminus L)\cap[\sqrt{10}, \sqrt{12}]\right)$, our claim about the dimension of $(M\setminus L)\cap [\sqrt{5}, \sqrt{12}])$ is reduced to proving\footnote{Contrary to the crude bound on the dimension of $(M\setminus L)\cap[3,\sqrt{10}]$ via the dimension of $M\cap[3,\sqrt{10}]$, the estimate $\textrm{dim}((M\setminus L)\cap[\sqrt{10}, \sqrt{12}]) < 0.93$ is non-trivial because, as we already mentioned, $
\textrm{dim}(M\cap[\sqrt{10}, \sqrt{12}])=1$ thanks to the fact that $\textrm{dim}(C(2))>1/2$.} that 
\begin{equation}\label{e.MM-upper-bound}
\textrm{dim}((M\setminus L)\cap[\sqrt{10}, \sqrt{12}]) < 0.93
\end{equation}

The crucial idea at this point is to put restrictions on the past or future shift dynamics of sequences $\theta\in\{1,2\}^{\mathbb{Z}}$ leading to Markov values in $(M\setminus L)\cap[\sqrt{10}, \sqrt{12}]$. For this sake, we take inspiration from the so-called \emph{shadowing lemma}, a classical result from the theory of uniformly hyperbolic dynamical systems asserting that pseudo-orbits are tracked by genuine orbits (cf. Hasselblatt--Katok book \cite{MR1326374}). In our setting, this lemma essentially says that if $m(\theta)\in (M\setminus L)\cap[\sqrt{10}, \sqrt{12}]$, then, up to transposition, the future dynamics of $\theta$ lives in the gaps of the \emph{symmetric block} $B=\{11,22\}^{\mathbb{Z}}\subset \{x\in\{1,2\}^{\mathbb{Z}}: x \textrm{ not containing }121\}$. Let us now try to make this statement more precise. 

We have a setting which can also be used in the analysis of other portions of $M\setminus L$: more concretely, we dispose of a symmetric block $C$, which here is $\{1,2\}^{\mathbb{Z}}$, together with its corresponding Gauss--Cantor set $K(C)=C(2)$, and a smaller symmetric block $B$ giving rise to the Gauss--Cantor set $K(B)=K(\{11,22\})$. The elements of $(M\setminus L)\cap[\sqrt{10}, \sqrt{12}]$ will be of the form $f(\theta)$ for some $\theta\in C$ such that $m(\theta)=f(\theta)$, i.e., let us consider 
$$Y=\{\theta\in C: m(\theta)=f(\theta)\in (M\setminus L)\cap[\sqrt{10}, \sqrt{12}]\}.$$ In order to prove that $\textrm{dim}((M\setminus L)\cap[\sqrt{10}, \sqrt{12}])\leq d$, it is enough (by taking countable subcoverings) to prove that, for every $\theta=(\theta_j)_{j\in \mathbb Z} \in Y$, there is $N\in \mathbb N$ such that 
$$\textrm{dim}(f(V_N(\theta)\cap Y))\le d,$$ where $V_N(\theta)=\{\tilde\theta=(\tilde\theta_j)_{j\in \mathbb Z}; \tilde\theta_j=\theta_j, -N\le j\le N\}$.

We say that a finite sequence $\omega=(\theta_j)_{-N\le j\le n}$ has an {\it allowed continuation} $\tau$, where $\tau$ is a finite string of $1$ and $2$, if there are sequences $\alpha, \beta\in \{1,2\}^{\mathbb N}$ such that $\beta^T \omega \tau \alpha\in Y$. 

We \emph{claim} that a shadowing lemma type argument implies that, if $\theta\in Y$, then (perhaps replacing $\theta$ by its transpose $\theta^T$) there is $N\in\mathbb N$ such that, for every $n\ge N$, the bifurcation tree (of allowed continuations) of such a $\theta$ is severely constrained: 

\begin{itemize}
\item either $\theta_{-N}\dots\theta_n$ has an unique 1-element continuation $\theta_{n+1}$,
\item or its (3-elements) allowed continuations are $112$ and $221$. 
\end{itemize}

In order to show this, observe that if $\theta\in Y$, then, since $L$ is closed, there is $N\in\mathbb N$ such that the distance of $m(\theta)$ to $L$ is larger than $\frac1{2^{N-4}}$. Also, if $\tilde\theta\in V_N(\theta)\cap Y$, then $|m(\tilde\theta)-m(\theta)|=|f(\tilde\theta)-f(\theta)|<\frac1{2^{N-2}}$ and $\lambda_j(\tilde\theta):=f(\sigma^j(\tilde\theta))\le f(\tilde\theta),\forall j\in\mathbb Z$. 

Suppose that the claimed statement above is not true for this $N$. Then there are allowed continuations beginning by $1$ and $2$, and we may assume without loss of generality that there is $\tau=(1,a,b)\ne (1,1,2)$, (the case of $(2,2,1)$ is analogous), $\alpha,\beta\in \{1,2\}^{\mathbb N}$ such that $\beta^T \omega \tau \alpha\in Y$, where $\omega=(\theta_j)_{-N\le j\le n}$. Since the smallest continued fractions in $C(2)$ beginning by $[0;1]$ begin indeed by $[0;1,1,2]$, and there are elements of $K(B)$ beginning by $[0;1,1,2]$, it follows that there is $\tilde\alpha\in\{11,22\}^{\mathbb N}$ such that 
$$[0;\tilde\alpha]\in K(B)\cap ([0;2,\overline{2,1}],[0;\tau \alpha]).$$ 
By monotonicity of maps of the type $h(y)=[c_0,c_1,\dots,c_k+y]$ (and noticing that $[0;2,\overline{2,1}]$ is the largest continued fraction in $C(2)$ beginning by $[0;2]$), it follows that, for $-N\le j\le n+1$, 
$$\lambda_j(\beta^T \omega \tilde\alpha)<m(\theta)+\frac1{2^{N-2}}+\frac1{2^{N-1}}<m(\theta)+\frac1{2^{N-3}}.$$ On the other hand, since $\tilde\alpha\in\{11,22\}^{\mathbb N}$, we have, for $j>n+1$, $\lambda_j(\beta^T \omega \tau \tilde\alpha)\le [2;\overline{1,1}] + [0;2,\overline{2,1}] < 3.0407<\sqrt{10}\le m(\theta)$. 

Hence, if the claimed statement above on allowed continuations fails both for $\theta$ and $\theta^T$, we may find $n, n'\ge N$ and $\tilde\alpha, \hat\alpha\in\{11,22\}^{\mathbb N}$ such that $\lambda_j(\hat\alpha^T \hat\omega \tilde\alpha)<m(\theta)+\frac1{2^{N-3}}$ for all $j\in \mathbb Z$, where $\hat\omega=(\theta_j)_{-n'\le j\le n}$. 

This leads to a contradiction: if $\tilde\alpha=(\tilde a_1,\tilde a_2,\dots)$ and $\hat\alpha=(\hat a_1,\hat a_2,\dots)$, defining \hfill\break $\check \alpha=(\tilde a_1,\tilde a_2,\dots,\tilde a_{2N},\hat a_{2N},\dots,\hat a_2,\hat a_1)\in\{1,2\}^{4N}$, then, if $\gamma=\overline{\hat\omega\check \alpha}$ is the periodic sequence with period $\hat\omega\check \alpha$, we have $|m(\gamma)-m(\theta)|<\frac1{2^{N-4}}$ and $m(\gamma)\in L$, an absurdity.

The result above on allowed continuations can be used to get that the set $\{[\theta_0;\theta_1,\theta_2,\dots]\}$ related to those allowed continuations of $\theta_{-N}\dots\theta_N$ is contained in a small ``Cantor set of gaps" $K_G$ and, \emph{a fortiori}, we have $(M\setminus L)\cap[\sqrt{10}, \sqrt{12}]\subset C(2)+K_G$. 

The non-trivial control of the bifurcation tree allows us to infer that 
$$\textrm{dim}(K_G)\leq s_0$$  
where $s_0\in(0,1)$ is any number such that 
\begin{eqnarray*} 
|I(\theta_1,\dots,\theta_n,1,1,2)|^{s_0}&+&|I(\theta_1,\dots,\theta_n,2,2,1)|^{s_0} \\ &\leq& |I(\theta_1,\dots,\theta_n)|^{s_0}
\end{eqnarray*}
for all $(\theta_0,\dots,\theta_n)\in\{1,2\}^n$. On the other hand, it is possible to derive\footnote{Here, we are using two useful facts: the box dimension of a Gauss--Cantor set (such as $C(2)$) coincides with its Hausdorff dimension and $\textrm{dim}(X\times Y)\leq \textrm{dim}_{box}(X)+\textrm{dim}(Y)$ where $\textrm{dim}_{box}(X)$ is the box dimension of $X$.} from the inclusion $(M\setminus L)\cap[\sqrt{10}, \sqrt{12}]\subset C(2)+K_G$ that 
$$\textrm{dim}((M\setminus L)\cap[\sqrt{10}, \sqrt{12}])\leq \textrm{dim}(C(2))+\textrm{dim}(K_G).$$

These facts permit to give a good upper bound on $\textrm{dim}((M\setminus L)\cap[\sqrt{10}, \sqrt{12}])$ because  $\textrm{dim}(C(2))<0.531281$ and some elementary estimates on continued fractions yield that $s_0=0.174813$ satisfies the above requirement, so that 
\begin{eqnarray*}
\textrm{dim}((M\setminus L)\cap[\sqrt{10}, \sqrt{12}])&<& 0.531281 + 0.174813 \\ &=& 0.706094.
\end{eqnarray*}
This completes the proof of \eqref{e.MM-upper-bound}.

\section{Beyond the classical spectra} Partly inspired by Perron's characterization of the classical spectra, several authors (including Maucourant, Paulin, Parkkonen and the second author) proposed dynamical generalizations of the Markov and Lagrange spectra. In a nutshell, \emph{dynamical} Lagrange and Markov spectra are obtained after replacing $\sigma$ by a \emph{general} dynamical system and $f$ by a general height function: for instance, given a homeomorphism $\varphi:M\to M$ of a topological space $M$, a compact $\varphi$-invariant subset $\Lambda$ of $M$, and a continuous function $f:M\to\mathbb{R}$, the \emph{dynamical Lagrange and Markov spectra}\footnote{An analogous definition can be made when the discrete time dynamical system $\varphi:M\to M$ is replaced by a continuous time dynamical system $(\varphi^t)_{t\in\mathbb{R}}$.} associated to $(f,\Lambda)$ are  
$$L(f,\Lambda):=\left\{\limsup_{n\to\infty}f(\varphi^{n}(x)):x\in \Lambda\right\}$$ 
and  
$$M(f,\Lambda)=\left\{ \sup_{n\in \mathbb{Z}}f(\varphi^{n}(x)):x\in \Lambda\right\}.$$ 

A direct generalization of the classical spectra is provided by the dynamical spectra associated to a diffeomorphism $\varphi:M\to M$ of a surface $M$ acting on a \emph{horseshoe}\footnote{This is a compact, $\varphi$-invariant, uniformly hyperbolic set of saddle type: cf. Hasselblatt--Katok book \cite{MR1326374}.} and a typical differentiable real function $f$. In fact, arbitrarily large compact parts of the classical Markov and Lagrange spectra can be viewed as dynamical Markov and Lagrange spectra associated to horseshoes of \emph{conservative} (i.e., area-preserving) diffeomorphisms. More precisely, as it is explained in \cite{MR1279059}, for each $m\geq 2$, the map $T_1:(0,1)\times(0,1)\to [0,1)\times(0,1)$ given by 
$$T_1(x,y)=\left(\left\{\frac1{x}\right\},\frac1{y+\lfloor 1/x\rfloor}\right)$$ 
preserves a smooth area-form near the horseshoe $\Lambda(m)=C(m)\times C(m)$ corresponding to the maximal invariant set of $(\frac1{m+1},1)\times(0,1)$. In particular, since  
$$T_1([0;a_0,a_1,a_2,\dots],[0;b_1,b_2,b_3,\dots])=$$
$$=([0;a_1,a_2,a_3,\dots],[0;a_0,b_1,b_2,\dots]),$$
we see that $T_1$ is a (piecewise) smooth, conservative realization of the natural extension of the Gauss map (compare with \eqref{e.Gauss-shift}). The dynamical Markov and Lagrange spectra of $(T_1,\Lambda(m))$ with respect to the function $f(x,y)=y+\frac1{x}$ have the same intersections with $(-\infty,m+1]$ as the classical Markov and Lagrange spectra. 

In 2018, Cerqueira and the authors established the continuity of the Hausdorff dimension across the dynamical Lagrange and Markov spectra of typical \emph{thin} horseshoes of {\it conservative} surface diffeomorphisms with respect to typical smooth functions (in analogy with the main continuity result in \S\ref{ss.main-results} for $L$ and $M$). More precisely, let $\varphi_0$ be a smooth diffeomorphism of a surface $M^2$ preserving an area-form $\omega$. Suppose that $\varphi_0$ possesses a thin horseshoe $\Lambda_0$ in the sense that its Hausdorff dimension is $\textrm{dim}(\Lambda_0)<1$. Denote by $\mathcal{U}$ a small $C^{\infty}$ neighborhood of $\varphi_0$ in the space $\textrm{Diff}_{\omega}^{\infty}(M)$ of smooth area-preserving diffeomorphisms of $M$ such that $\Lambda_0$ admits a \emph{continuation}\footnote{I.e., if $U_0$ is a neighborhood of $\Lambda_0$ such that $\Lambda_0=\bigcap\limits_{n\in\mathbb{Z}}\varphi_0^n(U_0)$, then $\mathcal{U}$ is taken small enough so that $\Lambda=\bigcap\limits_{n\in\mathbb{Z}}\varphi^n(U_0)$ still is a horseshoe for any $\varphi\in\mathcal{U}$.} $\Lambda$ for every $\varphi\in\mathcal{U}$. If $\mathcal{U}\subset\textrm{Diff}_{\omega}^{\infty}(M)$ is sufficiently small, then there exists a Baire residual subset $\mathcal{U}^{**}\subset \mathcal{U}$ with the following property. For every $\varphi\in\mathcal{U}^{**}$ and $r\geq 2$, there exists a $C^r$-open and dense subset $\mathcal{R}_{\varphi,\Lambda}\subset C^r(M,\mathbb{R})$ such that $$\textrm{dim}(L(\Lambda, f)\cap (-\infty, t)) = \textrm{dim}(M(\Lambda, f)\cap (-\infty, t))$$ 
is a continuous function of $t$ whenever $f\in \mathcal{R}_{\varphi,\Lambda}$.

Still concerning the beginning of the dynamical spectra, the second author \cite{MR4142446} proved that, for typical pairs $(f,\Lambda)$ as above, the minima of the corresponding Lagrange and Markov dynamical spectra coincide and are given by the image of a periodic point of the dynamics by the real function, solving a question by Yoccoz. 

Recently, Lima and the second author\footnote{See the article \emph{Phase transitions on the Markov and Lagrange dynamical spectra}, to appear in  Ann. Inst. H. Poincar\'e Anal. Non Lin\'eaire.} proved that, for typical pairs $(f,\Lambda)$ as in the previous paragraphs, 
$$\sup\{t\in\mathbb{R}: \textrm{dim}(M(f,\Lambda)\cap(-\infty,t))<1\}=$$
$$=\inf\{t\in\mathbb{R}: \textrm{int}(L(f,\Lambda))\cap(-\infty,t)\ne\emptyset\},$$
and, inspired by this result, they conjectured that the classical Lagrange spectrum must have non-empty interior right after the transition point where the classical Markov spectrum acquires Hausdorff dimension one, i.e., 
$$\textrm{int}(L)\cap (t_1,t_1+\epsilon) \neq\emptyset$$ for all $\epsilon>0$, where $t_1=\inf\{t: \textrm{dim}(M\cap[\sqrt{5},t])=1\}$.

Concerning the ending of dynamical spectra associated to horseshoes, Roma\~na and the second author proved that if $\Lambda$ is a (not necessarily conservative) horseshoe associated to a $C^2$-diffeomorphism $\varphi$ such that  $\textrm{dim}(\Lambda)>1$, then there is, arbitrarily close to $\varphi$ a diffeomorphism $\varphi_{0}$ and a $C^{2}$-neighborhood $W$ of $\varphi_{0}$ such that, if $\Lambda_{\psi}$ denotes the continuation of $\Lambda$ associated to $\psi\in W$, there is an open and dense set $H_{\psi}\subset C^{1}(M,\mathbb R)$ such that for all $f\in H_{\psi}$, we have 
$$
\textrm{int } L(f,\Lambda_{\psi})\neq\emptyset \textrm{ and } \textrm{int }M(f,\Lambda_{\psi})\neq\emptyset.
$$

Another direct generalization of the classical spectra for geodesic flows on negatively curved manifolds and moduli spaces of translation surfaces (along the lines of \S\ref{ss.modular-surface}) was studied by Maucourant, Paulin, Parkkonen, Artigiani, Delecroix, Hubert, Leli\`evre, Marchese, and Ulcigrai (among others): in their respective settings, these authors showed that their spectra shared some properties of the classical spectra such as isolated minima and a Hall's ray. We refer the reader to our recent book with Lima and Roma\~na \cite{LMMRn20} (and the references therein) for more details about dynamical generalizations of the classical spectra. 


\end{document}